\documentclass[12pt]{amsart}
\usepackage{amssymb}

\begin{document}

\def\dbl{[\hskip -1pt[}
\def\dbr{]\hskip -1pt]}
\title[Lie group structures on automorphism groups of CR manifolds]{Lie group structures on automorphism
groups of real-analytic CR manifolds}
\author{Bernhard Lamel}
\address{Universit\"at Wien, Fakult\"at f\"ur Mathematik, Nordbergstrasse 15, A-1090 Wien, \"Osterreich}
\email{lamelb@member.ams.org}%
\thanks{The first author was supported by the FWF, Projekt P17111.}
\author{Nordine Mir}
\address{Universit\'e de Rouen, Laboratoire de Math\'ematiques Rapha\"el Salem, UMR 6085 CNRS, Avenue de
l'Universit\'e, B.P. 12, 76801 Saint Etienne du Rouvray, France}
\email{Nordine.Mir@univ-rouen.fr}
\author{Dmitri Zaitsev}
\address{ School of Mathematics, Trinity College, Dublin 2,
Ireland} \email{zaitsev@maths.tcd.ie} 
\subjclass[2000]{32H02, 32V20, 32V25, 32V35, 32V40, 22E99} \keywords{Global CR automorphism group,
Lie group structure, jet parametrization}


\def\Label#1{\label{#1}}
\def\1#1{\ov{#1}}
\def\2#1{\widetilde{#1}}
\def\6#1{\mathcal{#1}}
\def\4#1{\mathbb{#1}}
\def\3#1{\widehat{#1}}
\def\K{{\4K}}
\def\LL{{\4L}}

\def \MM{{\4M}}
\def \S{{\4S}^{2N'-1}}

\def \B{{\4B}^{2N'-1}}

\def \H{{\4H}^{2l-1}}

\def \F{{\4H}^{2N'-1}}

\def \LL{{\4L}}

\def\Re{{\sf Re}\,}
\def\Im{{\sf Im}\,}
\def\id{{\sf id}\,}

\def\s{s}
\def\k{\kappa}
\def\ov{\overline}
\def\span{\text{\rm span}}
\def\ad{\text{\rm ad }}
\def\tr{\text{\rm tr}}
\def\xo {{x_0}}
\def\Rk{\text{\rm Rk\,}}
\def\sg{\sigma}
\def \emxy{E_{(M,M')}(X,Y)}
\def \semxy{\scrE_{(M,M')}(X,Y)}
\def \jkxy {J^k(X,Y)}
\def \gkxy {G^k(X,Y)}
\def \exy {E(X,Y)}
\def \sexy{\scrE(X,Y)}
\def \hn {holomorphically nondegenerate}
\def\hyp{hypersurface}
\def\prt#1{{\partial \over\partial #1}}
\def\det{{\text{\rm det}}}
\def\wob{{w\over B(z)}}
\def\co{\chi_1}
\def\po{p_0}
\def\fb {\bar f}
\def\gb {\bar g}
\def\Fb {\ov F}
\def\Gb {\ov G}
\def\Hb {\ov H}
\def\zb {\bar z}
\def\wb {\bar w}
\def \qb {\bar Q}
\def \t {\tau}
\def\z{\chi}
\def\w{\tau}
\def\Z{\zeta}
\def\phi{\varphi}
\def\eps{\epsilon}

\def \T {\theta}
\def \Th {\Theta}
\def \L {\Lambda}
\def\b {\beta}
\def\a {\alpha}
\def\o {\omega}
\def\l {\lambda}

\def \im{\text{\rm Im }}
\def \re{\text{\rm Re }}
\def \Char{\text{\rm Char }}
\def \supp{\text{\rm supp }}
\def \codim{\text{\rm codim }}
\def \Ht{\text{\rm ht }}
\def \Dt{\text{\rm dt }}
\def \hO{\widehat{\mathcal O}}
\def \cl{\text{\rm cl }}
\def \bS{\mathbb S}
\def \bK{\mathbb K}
\def \bD{\mathbb D}
\def \bC{\mathbb C}
\def \bL{\mathbb L}
\def \bZ{\mathbb Z}
\def \bN{\mathbb N}
\def \scrF{\mathcal F}
\def \scrK{\mathcal K}
\def \mc #1 {\mathcal {#1}}
\def \scrM{\mathcal M}
\def \cR{\mathcal R}
\def \scrJ{\mathcal J}
\def \scrA{\mathcal A}
\def \scrO{\mathcal O}
\def \scrV{\mathcal V}
\def \scrL{\mathcal L}
\def \scrE{\mathcal E}
\def \hol{\text{\rm hol}}
\def \aut{\text{\rm aut}}
\def \Aut{\text{\rm Aut}}
\def \J{\text{\rm Jac}}
\def\jet#1#2{J^{#1}_{#2}}
\def\gp#1{G^{#1}}
\def\gpo{\gp {2k_0}_0}
\def\emmp {\scrF(M,p;M',p')}
\def\rk{\text{\rm rk\,}}
\def\Orb{\text{\rm Orb\,}}
\def\Exp{\text{\rm Exp\,}}
\def\Span{\text{\rm span\,}}
\def\d{\partial}
\def\D{\3J}
\def\pr{{\rm pr}}

\def \CZZ {\C \dbl Z,\zeta \dbr}
\def \D{\text{\rm Der}\,}
\def \Rk{\text{\rm Rk}\,}
\def \CR{\text{\rm CR}}
\def \ima{\text{\rm im}\,}
\def \I {\mathcal I}

\def \M {\mathcal M}

\newtheorem{Thm}{Theorem}[section]
\newtheorem{Cor}[Thm]{Corollary}
\newtheorem{Pro}[Thm]{Proposition}
\newtheorem{Lem}[Thm]{Lemma}

\theoremstyle{definition}\newtheorem{Def}[Thm]{Definition}

\theoremstyle{remark}
\newtheorem{Rem}[Thm]{Remark}
\newtheorem{Exa}[Thm]{Example}
\newtheorem{Exs}[Thm]{Examples}

\numberwithin{equation}{subsection}

\def\bl{\begin{Lem}}
\def\el{\end{Lem}}
\def\bp{\begin{Pro}}
\def\ep{\end{Pro}}
\def\bt{\begin{Thm}}
\def\et{\end{Thm}}
\def\bc{\begin{Cor}}
\def\ec{\end{Cor}}
\def\bd{\begin{Def}}
\def\ed{\end{Def}}
\def\be{\begin{Exa}}
\def\ee{\end{Exa}}
\def\bpf{\begin{proof}}
\def\epf{\end{proof}}
\def\ben{\begin{enumerate}}
\def\een{\end{enumerate}}

\newcommand{\dbar}{\bar\partial}
\newcommand{\genmat}{\lambda}
\newcommand{\polynorm}[1]{{|| #1 ||}}
\newcommand{\vnorm}[1]{\left\|  #1 \right\|}
\newcommand{\asspol}[1]{{\mathbf{#1}}}
\newcommand{\Cn}{\mathbb{C}^n}
\newcommand{\Cd}{\mathbb{C}^d}
\newcommand{\Cm}{\mathbb{C}^m}
\newcommand{\C}{\mathbb{C}}
\newcommand{\CN}{\mathbb{C}^N}
\newcommand{\CNp}{\mathbb{C}^{N^\prime}}
\newcommand{\Rd}{\mathbb{R}^d}
\newcommand{\Rn}{\mathbb{R}^n}
\newcommand{\RN}{\mathbb{R}^N}
\newcommand{\R}{\mathbb{R}}
\newcommand{\bR}{\mathbb{R}}
\newcommand{\N}{\mathbb{N}}
\newcommand{\dop}[1]{\frac{\partial}{\partial #1}}
\newcommand{\vardop}[3]{\frac{\partial^{|#3|} #1}{\partial {#2}^{#3}}}
\newcommand{\br}[1]{\langle#1 \rangle}
\newcommand{\infnorm}[1]{{\left\| #1 \right\|}_{\infty}}
\newcommand{\onenorm}[1]{{\left\| #1 \right\|}_{1}}
\newcommand{\deltanorm}[1]{{\left\| #1 \right\|}_{\Delta}}
\newcommand{\omeganorm}[1]{{\left\| #1 \right\|}_{\Omega}}
\newcommand{\nequiv}{{\equiv \!\!\!\!\!\!  / \,\,}}
\newcommand{\bk}{\mathbf{K}}
\newcommand{\p}{\prime}
\newcommand{\tV}{\mathcal{V}}
\newcommand{\poly}{\mathcal{P}}
\newcommand{\ring}{\mathcal{A}}
\newcommand{\ringk}{\ring_k}
\newcommand{\ringktwo}{\mathcal{B}_\mu}
\newcommand{\germs}{\mathcal{O}}
\newcommand{\On}{\germs_n}
\newcommand{\mcl}{\mathcal{C}}
\newcommand{\formals}{\mathcal{F}}
\newcommand{\Fn}{\formals_n}
\newcommand{\autM}{{\Aut (M,0)}}
\newcommand{\autMp}{{\Aut (M,p)}}
\newcommand{\holmaps}{\mathcal{H}}
\newcommand{\biholmaps}{\mathcal{B}}
\newcommand{\autmaps}{\mathcal{A}(\CN,0)}
\newcommand{\jetsp}[2]{ G_{#1}^{#2} }
\newcommand{\njetsp}[2]{J_{#1}^{#2} }
\newcommand{\jetm}[2]{ j_{#1}^{#2} }
\newcommand{\glnc}{\mathsf{GL_n}(\C)}
\newcommand{\glmc}{\mathsf{GL_m}(\C)}
\newcommand{\glc}{\mathsf{GL_{(m+1)n}}(\C)}
\newcommand{\glk}{\mathsf{GL_{k}}(\C)}
\newcommand{\smC}{\mathcal{C}^{\infty}}
\newcommand{\anC}{\mathcal{C}^{\omega}}
\newcommand{\kC}{\mathcal{C}^{k}}



\begin{abstract}
Given any real-analytic CR manifold $M$, we
provide general conditions on $M$ guaranteeing that the group of
all its global real-analytic CR automorphisms $\Aut_\CR(M)$ is a Lie
group (in an appropriate topology). 
In particular, we obtain a Lie group structure for
$\Aut_\CR(M)$ when $M$ is an arbitrary compact real-analytic
hypersurface embedded in some Stein manifold.
\end{abstract}

\maketitle

\section{Introduction}\Label{int}

There exist a wide variety of results concerned with the structure
of the automorphism group of a given geometric structure. In
Riemannian Geometry, the classical Myers-Steenrod theorem
\cite{MS39} states that the group of all isometries of a Riemannian
manifold is a Lie group. 
 H.\ Cartan \cite{HCa} proved an analogous result for the group of holomorphic
automorphisms of a bounded domain in $\C^N$. 
Cartan's techniques have in turn been used to establish general
results for groups of diffeomorphisms of real or complex manifolds,
see e.g.\ \cite{BM0}.

 In this paper, we consider an analogous question for CR
manifolds (that one can think of as a boundary or CR version of Cartan's
Theorem mentioned above): 

{\em Under what conditions on a real-analytic CR manifold $M$ is
the group $\Aut_\CR (M)$ of all real-analytic CR automorphisms of $M$
a Lie group in an appropriate topology?}

Here for every $r\in \N \cup \{\infty,\omega\}$, we equip $\Aut_\CR (M)$ 
with a natural topology that we call  ``compact-open $\6C^r$ topology'', which is defined as follows.
For open subsets $\Omega \subset \R^n$ and $\Omega' \subset \R^{n'}$,
consider the space $\6C^r(\Omega,\Omega')$ of all maps of class $\6C^r$ from $\Omega$ to $\Omega'$.
If $r\in \N \cup \{\infty\}$, $\6C^r(\Omega,\Omega')$ is equipped with the topology of uniform convergence
on compacta together with all partial derivatives of order up to $r$.
In case $r=\omega$, the space $\anC(\Omega,\Omega')$
is equipped with its topology as an inductive limit
of Fr\'echet spaces of holomorphic maps between open neighborhoods
of $\Omega$ and $\Omega'$ in $\C^n$ and $\C^{n'}$ respectively.
The compact-open  $\6C^r$ topology on $\Aut_\CR (M)$ is now induced by 
the appropriate topology relative to the coordinate charts for the maps and their inverses
(see e.g.\ \cite{BRWZ} for a more detailed discussion).
For brevity, we adopt the order $k<\infty<\omega$ for any integer $k$.

In this paper, we exhibit general
sufficient conditions on $M$ that provide an affirmative answer to the above question. 
We begin with the following special case of our more general results,
that is particularly easy to state:

\begin{Cor}\label{t:lieh}
Let $M$ be a compact real-analytic hypersurface in a Stein manifold
of complex dimension at least two. 
Then the group $\Aut_\CR(M)$ of
all $($global$)$ real-analytic CR automorphisms of $M$ is a Lie
group in the compact-open $\anC$ topology and the
action $\Aut_\CR(M)\times M \to M$ is real-analytic.
Furthermore, the compact-open $\6C^r$ topologies on $\Aut_\CR(M)$
coincide for 
$r = \infty,\omega$ and 
$r\ge k$, where $k$ is an integer
depending only on $M$.
\end{Cor}

Corollary~\ref{t:lieh} is a direct consequence of the following more
general result that also applies to CR manifolds of higher codimension.
In the following statement, the notions of essential finiteness, finite
nondegeneracy and minimality must be understood in the sense of \cite{BERbook} 
(see also Section~\ref{s:perturb} for more details).

\begin{Thm}\label{t:unify}
Let $M$ be a real-analytic CR manifold. Assume that $M$
has finitely many connected components, is minimal everywhere and
that there exists a compact subset $K\subset M$ such that:
\begin{enumerate}
\item[(i)] $M$ is essentially finite at all points of $K$;
\item[(ii)] $M$ is finitely nondegenerate at all points of
$M\setminus K$.
\end{enumerate}
Then $\Aut_\CR(M)$ is a Lie group in the compact-open $\anC$
 topology and the
action $\Aut_\CR(M)\times M \to M$ is real-analytic.
Furthermore, the compact-open $\6C^r$ topologies on $\Aut_\CR(M)$
coincide for 
$r = \infty,\omega$ and $r\ge k$ for some integer $k$,
where $k$ is an integer
depending only on $M$.
\end{Thm}

Theorem~\ref{t:unify} provides a generalization of all known
corresponding results for real-analytic CR manifolds. It also covers
new situations, such as in Corollary~\ref{t:lieh}; indeed, any
real-analytic compact hypersurface in a Stein manifold is
essentially finite and minimal at {\em each} of its points (see
e.g.\ \cite{DF, BERbook}).

For the case of real hypersurfaces whose
Levi form is nondegenerate at every point, the conclusion of
Theorem~\ref{t:unify} follows from the work of E.\ Cartan
\cite{Ca1, Ca2}, Chern-Moser \cite{CM}, Tanaka \cite{Ta} and
Burns-Schnider \cite{BS}. For the
case of Levi-degenerate CR manifolds, the same conclusion was
recently obtained by Baouendi, Rothschild, Winkelmann and the third author  
\cite{BRWZ} for the class of finitely
 nondegenerate minimal CR manifolds, which corresponds here to our
 Theorem~\ref{t:unify} with $K=\emptyset$. (We should point out that
 the results in those mentioned papers also apply for merely smooth CR manifolds as well, 
based on the previous work \cite{KZ05},
but in this paper we shall focus on the real-analytic category.)

 In addition to the compact hypersurface case considered in Corollary~\ref{t:lieh},
 an important class of CR manifolds for which the previously
 known results do not apply and for which the conclusion of
Theorem~\ref{t:unify} holds
 is that of {\em compact} minimal real-analytic CR submanifolds
 embedded in a Stein manifold.
 Again, the condition of essential finiteness holds here at every point, see \cite{DF, BERbook}
(whereas the condition of finite nondegeneracy
 holds only outside a proper real-analytic subvariety which need not be empty in general); 
 Therefore taking $K=M$ in Theorem~\ref{t:unify}, we
 obtain the following extension of Corollary~\ref{t:lieh} to higher codimension:

\begin{Cor}\label{t:liehighcodim}
Let $M$ be a compact real-analytic CR submanifold in a Stein manifold. Assume that $M$ is
minimal at every point. Then the group $\Aut_\CR(M)$ of all
$($global$)$ real-analytic CR automorphisms of $M$ is a Lie group in
the compact-open $\anC$ topology and the action
$\Aut_\CR(M)\times M \to M$ is real-analytic.
Furthermore, the compact-open $\6C^r$ topologies on $\Aut_\CR(M)$
coincide for 
$r = \infty,\omega$ and 
$r\ge k$, where $k$ is an integer
depending only on $M$.
\end{Cor}

On the other hand, Theorem~\ref{t:unify}
also applies to cases with $M$ noncompact also not covered by previously known results. 
Let us illustrate this with an example:
\begin{Exa}
 The hypersurface $M\subset \C^2$ given by
$$|z|^2-|w|^4 = 1$$
is Levi-nondegenerate at all its points except the circle $S^1\times \{0\}\subset M$,
where $M$ is essentially finite.
Hence, Theorem~\ref{t:unify} applied with $K:= S^1\times \{0\}$,
yields that $\Aut_\CR (M)$ is a Lie group. 
On the other hand, $M$ is not finitely nondegenerate at any point of $K$
and hence the results from \cite{BRWZ} do not apply to $M$.
\end{Exa}

Our proof of Theorem~\ref{t:unify} makes use of the recent
developments providing a relationship between
 various notions and results concerning jet parametrization of local CR
diffeomorphisms \cite{BER97, Z97, BER99b, E, KZ05, LM05} and Lie
group structures on (local and) global groups of automorphisms of CR
manifolds \cite{BRWZ}. In the next section, we give in
Theorem~\ref{main-technical} new sufficient conditions on a
connected real-analytic CR manifold $M$, in terms of local jet
parametrization properties of CR automorphisms, that ensure that $\Aut_\CR(M)$
is a Lie group. Then the remainder of the paper is
devoted to prove that under the assumptions of Theorem~\ref{t:unify}
the conditions of Theorem~\ref{main-technical} are fulfilled. To
this end, we establish, following the analysis of the first two
authors' paper \cite{LM05}, a new parametrization theorem
(Theorem~\ref{t:jetparamthm}) for local CR automorphisms that may be
of independent interest. The proof of Theorem~\ref{t:unify} is given
in Section \ref{ss:final}. We conclude the paper by giving in
Section \ref{s:lmz} an alternative proof of
Corollary~\ref{t:liehighcodim} following \cite{Z97} that does not make use of
Theorem~\ref{main-technical} but requires compactness of the manifold $M$.

\bigskip

\noindent  {\bf Acknowledgement.} The third author would like to thank A. Isaev for helpful discussions.

\section{New sufficient conditions for the automorphism group being a Lie group}\label{s:jetparam}

Let $M$ be a real-analytic manifold and $k$ a positive integer. We
use the notation $G^k(M)$ for the fiber bundle of all $k$-jets of
local real-analytic diffeomorphisms of $M$. For every point $p\in M$, we denote by $G_p^k(M)$ the fiber of
$G^k(M)$ at $p$. Given a germ of a local real-analytic
diffeomorphism $h\colon (M,p)\to M$, we write $j^k_p h\in G_p^k(M)$
for the corresponding $k$-jet. For instance, $j^k_p\id$ is the
$k$-jet of the identity map of $M$, regarded as a germ at $p$. In
local coordinates, $j^k_p h$ is given by the source $p$, the target
$h(p)$ and the collection of all partial derivatives of $h$ at $p$
up to order $k$. (See e.g.\ \cite{GG} for more details on this terminology.)

We now fix an arbitrary set $\6S$ of germs of local real-analytic
diffeomorphisms $h\colon (M,p)\to M$ with possibly varying reference
point $p\in M$ and, as in \cite{BRWZ},  consider the following
condition.

\bd\Label{abstract-param} Let $k$ be a positive integer and $p_0\in
M$. We say that $\6S$ {\em has the real-analytic jet parametrization
property of order $k$ at $p_0$} if there exist open neighbourhoods
$\Omega'$ of $p_0$ in $M$, $\Omega''$ of $j^k_{p_0}{\sf id}$ in
$G^k(M)$ and a real-analytic map $\Psi \colon \Omega'\times
\Omega''\to M$ such that, for every germ $h\colon (M,p)\to M$ in
$\6S$ with $p\in \Omega'$ and $j_p^k h\in \Omega''$, the identity
$h(\cdot) \equiv \Psi (\cdot,j_p^k h)$ holds in the sense of germs
at $p$. \ed

The following theorem is one of the key ingredients in the proof of
Theorem~\ref{t:unify}.

\bt\Label{main-technical} Let $M$ be a connected real-analytic CR
manifold. Assume that there exist an integer $k$ and a compact
subset $K\subset M$ such that the following holds:
\begin{enumerate}
\item[(i)] For every $p_0\in K$, the set of all germs at $p_0$ of local CR diffeomorphisms of $M$
has the real-analytic jet parametrization property of order $k$ at
$p_0$;
\item[(ii)] The set of all germs at all points of local CR diffeomorphisms
has the real-analytic jet parametrization property at every
point $p_0\in M\setminus K$ of some finite order possibly depending on $p_0$.
\end{enumerate}
Then $\Aut_\CR(M)$ is a Lie group in the compact-open $\anC$
topology and the
action $\Aut_\CR(M)\times M \to M$ is real-analytic.
Furthermore, the compact-open $\6C^r$ topologies on $\Aut_\CR(M)$
coincide for 
$r = \infty,\omega$ and 
$r\ge k$, where $k$ is an integer
depending only on $M$.
\et

In the case $K=\emptyset$, Theorem~\ref{main-technical} is
contained in \cite{BRWZ}. Heuristically speaking the points of
$M\setminus K$ in Theorem~\ref{main-technical} (ii) fulfil a
"strong'' jet parametrization property (namely, a so-called complete
system in the sense of \cite{KZ05,BRWZ}). In
Theorem~\ref{main-technical}, we allow some points to satisfy a
weaker property (namely condition (i)), but we have to pay the price
by requiring that all these points lie in a compact subset of $M$.

\begin{proof}[Proof of Theorem~{\rm \ref{main-technical}}]
Let $K\subset M$ be the compact subset as in
Theorem~\ref{main-technical}.
We first apply the parametrization property for the set of all germs
at a fixed point $p_0\in K$, which holds in view of (i); without
loss of generality we may assume that $K$ is nonempty. By
Definition~\ref{abstract-param}, for every fixed $p_0\in K$, we can
find open neighbourhoods $\Omega'$ of $p_0$ in $M$, $\Omega''$ of
$j^k_{p_0}{\sf id}$ in $G^k(M)$ and a real-analytic map $\Psi \colon
\Omega'\times \Omega''\to M$ such that, for every germ $h\colon
(M,p_0)\to M$ of a local CR diffeomorphism of $M$ with $j_{p_0}^k h\in
\Omega''$, we have the identity $h(\cdot) \equiv \Psi
(\cdot,j_{p_0}^k h)$ in the sense of germs at $p_0$. Let $\2\Omega'$
(resp.\ $\2\Omega''$) be a smaller neighbourhood of $p_0$ in
$\Omega'$  (resp.\ of $j^k_{p_0}{\sf id}$ in $\Omega''$) which is
relatively compact in $\Omega'$ (resp.\ $\Omega''$), chosen for
every $p_0\in K$. Without loss of generality, all neighbourhoods here are
connected. Using the compactness of $K$ and passing to a finite
subcovering, we obtain a finite collection of points
$p_1,\ldots,p_s\in K$, the corresponding neighbourhoods
$$\Omega'_m \supset \supset \2\Omega'_m \ni p_m, \quad
\Omega''_m \supset \supset \2\Omega''_m \ni j^k_{p_m}{\sf id},$$ 
and real-analytic maps $\Psi_m \colon \Omega'_m\times \Omega''_m\to M$
for $m=1,\ldots,s$, such that $(\2\Omega'_m)$ is a covering of $K$.

We next define neighbourhoods $\6U$ and $\2{\6U}$ of the identity
mapping in $\Aut_\CR(M)$ with respect to the compact-open $\6C^k$ topology as follows:
\begin{equation}\Label{u-def}
\2{\6U} := \{ g\in \Aut_\CR(M) : j^k_{p_m} g \in \2\Omega''_m, \,
1\le m\le s\}, \quad {\6U} := \{ g\in \2{\6U} : g^{-1} \in \ \2{\6U}
\}.
\end{equation}
It is clear from the definition of the topology chosen
that $\6U$ is indeed an open set. 
Obviously the
same conclusion holds for the compact-open $\smC$ topology as well
as for the compact-open $\mathcal{C}^r$ topology for any
$r \ge k$.

Our main step of the proof will be to show that $\6U$ is relatively
compact in $\Aut_\CR(M)$. We shall prove it with respect to the
compact-open $\6C^k$ topology, which is Fr\'echet and hence, in particular,
metrizable. Thus it suffices to prove that the closure of $\6U$ is
sequentially compact. Let $(f_n)$ be any sequence in $\6U$, for
which we shall prove that there exists a convergent subsequence. In
view of \eqref{u-def}, we have
$$j^k_{p_m} f_n \in \2\Omega''_m \subset\subset \Omega''_m \subset G^k(M), \quad m=1,\dots,s,$$
for every $n$. Hence, passing to a subsequence, we may assume that
$j^k_{p_m} f_n$ converges to some $\L_m\in \Omega''_m$ for each
$m=1,\dots,s$.

Following the strategy of \cite{BRWZ}, we denote by $\6O$ the open
set of points $q\in M$ with the property that $(f_n)$ converges in
the compact-open $\anC$ (and hence any $\6C^r$ with $r\ge k$)
topology in a neighbourhood $V$ of $q$ in $M$ to a map
$f\colon V\to M$ such that the Jacobian of $f$ at $q$ is nonzero. We want to show
that $\6O$ is nonempty and closed in $M$. By our construction, we
have $f_n(\cdot) \equiv \Psi_1 (\cdot,j_{p_1}^k f_n)$ in the sense
of germs at $p_1$ and hence, by the identity principle for
real-analytic functions, all over $\Omega'_1$. Since $j^k_{p_1} f_n$
converges to $\L_1 \in \Omega''_1$, it clearly follows that
$f_n|_{\Omega'_1}\to f:=\Psi_1(\cdot,\Lambda_1)$ as $n\to +\infty$
in the compact-open $\anC$ topology on $\Omega'_1$. In particular, we also have
$j^k_{p_1}f_n\to j_{p_1}^kf$ and since $j_{p_1}^kf \in
\Omega''_1\subset G^k(M)$, we immediately see that $p_1\in \6O$,
proving that $\6O$ is nonempty. To show that $\6O$ is closed,
let $q_0$ be any point in the closure of $\6O$ in $M$. We now
distinguish two cases.

{\bf Case 1}. $q_0\notin K$. Here we can repeat the arguments of the
proof of \cite[Lemma 3.3]{BRWZ} to show that $q_0\in \6O$.

{\bf Case 2}. $q_0\in K$. Here we only have the restricted
parametrization given by (i) and hence cannot use the same arguments
as in Case 1; instead, we use our construction. Since the
neighbourhoods $\2\Omega'_m$, $m=1,\ldots,s$, cover $K$, we have
$q_0\in \2\Omega'_{m_0} \subset \subset \Omega'_{m_0}$ for some
$m_0$ and let $p_{m_0}\in K$ be the corresponding point. The
sequence of the $k$-jets $\L_{m_0}^n:= j^k_{p_{m_0}} f_n$ converges
to $\L_{m_0}$ by our assumptions above
 and
therefore
\begin{equation}\Label{m0par}
  f_n(\cdot) \equiv \Psi_{m_0} (\cdot,\L_{m_0}^n) \to \Psi_{m_0} (\cdot,\L_{m_0}),
\end{equation}
which immediately implies that $q_0\in \6O$.

Summarizing, we have shown that $\6O$ is nonempty, open and closed
in $M$ and therefore $\6O=M$, i.e.\ $(f_n)$ converges on
$M$ to a real-analytic map $f\colon M\to M$ which is automatically CR. 
Furthermore, by our construction of $\6U$,
also the sequence of the inverses $f_n^{-1}$ is in $\6U$. Hence
similar arguments show that this sequence converges to another
real-analytic CR self-map $g$ of $M$. Then it follows that $g\circ f =
f\circ g = \id$ and therefore $f\in \Aut_\CR(M)$. This completes the
proof that the chosen neighbourhood  $\6U$ of $\id$ in $\Aut_\CR(M)$
is relatively compact. Since any  $g\in \Aut_\CR(M)$
has $g\6U$ as its neighbourhood, it follows that the whole group
$\Aut_\CR(M)$ is locally compact.

As in \cite{BRWZ}, we make use of the following theorem of
Bochner-Montgomery \cite{BM}, \cite[Theorem 2, p.~208]{MZ}: 
\bt[Bochner-Montgomery]\Label{bm1} Let
$G$ be a locally compact topological group acting effectively and continuously on a
smooth manifold $M$ by smooth diffeomorphisms. 
Then $G$ is a Lie group and the action
$G\times M\to M$ is smooth. 
\et

Indeed, we have just shown that $G:= \Aut_\CR(M)$ is locally
compact. 
Since the action $\Aut_\CR(M)\times M \to M$ is obviously effective, Theorem~\ref{bm1} shows that $\Aut_\CR(M)$ is
a Lie group and its action is smooth. 
The coincidence of the compact-open $\6C^r$ topologies on $\Aut_\CR(M)$ for $r\ge k$
also follows from the proof.
Finally the analyticity of the
action follows from another result of Bochner-Montgomery \cite{BM0}:

\bt[Bochner-Montgomery]\Label{bm2} Let $G$ be a Lie group acting
continuously on a real-analytic manifold $M$ by real-analytic
diffeomorphisms. Then the action $G\times M\to M$ is real-analytic.
\et

The proof of Theorem~\ref{main-technical} is complete.
\end{proof}

\section{Parametrization of local CR diffeomorphisms}\label{s:perturb}
In order to deduce Theorem~\ref{t:unify} from
Theorem~\ref{main-technical}, we will establish a jet
parametrization property of local CR diffeomorphisms for a certain
class of real-analytic CR submanifolds in complex space. Such a
property has already been established by the first two authors in
\cite{LM05} for an appropriate class of CR manifolds for local CR
diffeomorphisms {\em which furthermore fix a given point of the
manifold}. However, in view of Definition~\ref{abstract-param}, we
need to extend such a parametrization property to local CR
diffeomorphisms {\em which do not necessarily fix a base point}. In
what follows, we make the above statements precise and show how they
may be derived from the analysis given in the paper \cite{LM05}.

The class of germs of real-analytic generic submanifolds we shall
consider in this paper is the one introduced by the first two
authors in \cite{LM05}, denoted by ${\mathcal C}$, whose definition
we now recall. Denote by $(M,p_0)$ a germ of a real-analytic generic
submanifold of $\CN$ (or, more generally, of any complex manifold)
of CR dimension $n$ and real codimension $d$, i.e.\ $N=n+d$,
$T_{p_0}M + iT_{p_0}M = T_{p_0}\CN$ and $n=\dim_\C (T_{p_0}M\cap
iT_{p_0}M)$. Let $\rho =(\rho_1,\ldots,\rho_d)$ be a real-analytic
vector valued defining function for $M$ in some neighbourhood $U$ of
$p_0$ in $\CN$ satisfying $\partial \rho_1\wedge \ldots \wedge
\partial \rho_d\not =0$. Using standard notation, we write $\rho$ as
a convergent power series (after shrinking $U$ if necessary)
\[\rho(Z,\1{Z})=\sum_{\alpha,\beta\in\N^N} \rho_{\a\b}(Z-p_0)^\a
(\1{Z-p_0})^\b,\quad Z\in U,\] where $\rho_{\alpha,\beta}\in \C^d$
satisfy $\rho_{\alpha,\beta} = \overline{\rho_{\beta,\alpha}}$, and
complexify it to the power series
\[\rho(Z,\zeta)=\sum \rho_{\a\b}(Z-p_0)^\a (\zeta-\1p_0)^\b,
\quad \d\rho_1\wedge\cdots\wedge\d\rho_d\ne0, \] 
with
$(Z,\zeta)\in\CN\times\CN$, which we still denote by $\rho$. It is
easy to see that the complexification $\rho(Z,\zeta)$ is still
convergent in a suitable neighbourhood of $(p_0,\1p_0)$ that (after
shrinking $U$ again if necessary) can be chosen of the form
$U\times\1U\subset \CN\times \CN$. Recall that the {\em Segre
variety} $S_q$ of a point $q\in U$ is the $n$-dimensional complex
submanifold of $U$ given by $S_q:=\{Z\in U:\rho (Z,\1q)=0 \}$. 
Furthermore, the
{\em complexification of $M$} is defined to be the
$2n+d$-dimensional complex submanifold of $U\times\1U$ given by
\begin{equation}\label{complexification}
{\mathcal M}:=\{(Z,\zeta)\in U\times\1U :\rho (Z,\zeta)=0\} =
\left\{ (Z,\zeta)\in U\times\1U \colon Z \in S_{\bar \zeta}
\right\}.
\end{equation}
For every integer $k$ and for $q\in \CN$, we denote by
$J^{k,n}_q(\CN)$ the space of all jets at $q$ of order $k$ of
$n$-dimensional complex submanifolds of $\CN$ passing through $q$.
For every $q\in M$ sufficiently close to $p_0$, we consider the
anti-holomorphic map $\pi_q^k$ defined as follows:
\begin{equation}\label{e:segredef} \pi_q^k \colon S_q \to J^{k,n}_q(\CN), \quad \pi_q^k (\xi) = j^k_q
S_\xi,
\end{equation}
where $j^k_q S_\xi$ denotes the $k$-jet at $q$ of the submanifold
$S_\xi$ (see e.g.\ \cite{Z99} for more details on jets of
complex submanifolds used here, and also \cite{LM05}).

Following \cite{LM05}, we say that the germ $(M,p_0)$ belongs to the
class $\mcl$ if {\em the anti-holomorphic map $\pi^k_{p_0}$ is
generically of full rank $n= \dim S_{p_0}$ in any neighborhood of $p_0$, for $k$ sufficiently
large}. For $(M,p_0)\in \mcl$, we denote by $\kappa_M(p_0)$ the
smallest integer $k$ for which the map $\pi_{p_0}^k$ is of generic
rank $n$. Since the Segre varieties are associated to $(M,p_0)$ in a
biholomorphically invariant way, the integer $\kappa_M(p_0)$ is a
biholomorphic invariant of the germ $(M,p_0)$. Note
furthermore that the condition for a germ of real-analytic generic
submanifold to belong to the class $\mcl$ is an {\em open} condition
in the sense that, if $M\in \mcl$ is given by the equation
$\rho(Z,\bar Z)=0$ as above, then $\2M\in\mcl$ for any $\2M$ given
by the equation $\2\rho(Z,\bar Z)=0$ with $\2\rho$ sufficiently
close to $\rho$ in the $\smC$ topology (see e.g.\ \cite{GG} for
details on this topology; here it is enough to assume that
$\2\rho$ is close to $\rho$ in the $\kC$-topology for a suitable
$k$). In particular, there exists a neighbourhood $V$ of $p_0$ in
$M$ such that $(M,q)\in\mcl$ for all $q\in V$ and moreover, it is
clear from the definition that $\kappa_M(q)$ is {\em upper semi-continuous} on $V$.

We also recall that $M$ is {\em essentially finite} (resp.\
 {\em finitely nondegenerate}) at $p_0$ if the map $\pi^k_{p_0}$ is
 {\em finite} near $p_0$ (resp.\ an {\em immersion} at $p_0$) for
$k$ sufficiently large (see \cite{BHR, BERbook} for more details).
It follows  that finite nondegeneracy of $M$ at $p_0$ implies
essential finiteness of $M$ at $p_0$ which in turn implies that
$(M,p_0)\in \mcl$. Recall also that $M$ is {\em minimal} at $p_0$ if
there does not exist any CR submanifold of lower dimension contained in $M$ and
passing through $p_0$ with the same CR dimension as that of $M$ (see
\cite{Tu, BERbook}).

For  a real-analytic CR submanifold $M\subset \C^N$ that is not
necessarily generic and for a point $p_0\in M$, we say that
$(M,p_0)$ is in the class $\mathcal{C}$ if it is in the class
$\mathcal{C}$ when considered as a generic submanifold of its
intrinsic complexification, i.e.\ the minimal germ of a complex submanifold of $\C^N$
containing $(M,p_0)$
(see e.g.\ \cite{BERbook} for this
notion). Finally we should also note that the local nondegeneracy conditions defined above
are defined in the same way for abstract real-analytic CR manifolds since such manifolds can always be locally embedded  in some complex euclidean space $\C^q$ for some integer $q$, see e.g.\ \cite{BERbook}.

 Finally, we refer the reader to \cite{LM05} for examples
of manifolds that belong to the class $\mathcal{C}$, as well as for
a more thorough discussion of the relation between this
nondegeneracy condition and other well-known nondegeneracy
conditions such as essential finiteness and finite nondegeneracy. We
only stress in this paper the following fact that will be used
implicitly in the proofs of Corollaries~\ref{t:lieh} and \ref{t:liehighcodim} and that
follows from a result of \cite{DF} : {\em for every
compact real-analytic CR submanifold $\Sigma$ embedded in some Stein
manifold and for every $q\in \Sigma$, $\Sigma$ is essentially finite
at $q$, and in particular $(\Sigma, q)\in {\mathcal C}$} (see
\cite{LM05} for more details).

The following parametrization theorem is the second main ingredient
of the proof of Theorem~\ref{t:unify}.

\begin{Thm}\label{t:jetparamthm}
Let $M\subset \C^N$ be a real-analytic CR submanifold of codimension
$d$ and $p_0\in M$. Assume that $(M,p_0)$ is minimal and belongs to
the class ${\mathcal C}$
and set $\ell_0:=2(d+1)\kappa_M(p_0)$. Then the set of all germs
$h\colon (M,p_0)\to M$ of local CR diffeomorphisms of $M$  has the
real-analytic jet parametrization property of order $\ell_0$ at
$p_0$.
\end{Thm}

As mentioned above, the difference between Theorem
\ref{t:jetparamthm} and \cite[Theorem 7.3]{LM05} is due to the fact
that the parametrization theorem given in \cite{LM05} is obtained
for the set of germs $h\colon (M,p_0)\to (M,p_0)$ of local CR
diffeomorphisms with fixed source $p_0$ but also with fixed target
$p_0$. The version given here by Theorem~\ref{t:jetparamthm} allows
to parametrize local CR diffeomorphisms which send the point $p_0$ to
a varying target point $p\in M$ (close to $p_0$) that has to be
regarded as an additional parameter. We will provide a deformation
version of \cite[Theorem 7.3]{LM05} which allows us to treat this
additional parameter in Theorem~\ref{t:deformation} below.
(Note that it is not always possible to parametrize in a proper sense
the germs of all local CR diffeomorphisms with varying both source and targets, 
see Remark~\ref{cannot} below).

Before we proceed, we need to introduce some further
terminology. Given a real-analytic manifold $E$ and a point $p_0\in
\C^N$, a {\em real-analytic family of germs at $p_0$ of
real-analytic generic submanifolds} $(M_\eps)_{\eps\in E}$ of
$\C^N$, is given by a family of convergent power series mapping in $Z$ and $\bar Z$ centered at $p_0$,
$\rho(Z,\bar Z;\eps)= (\rho_1(Z,\bar Z;\eps),\ldots,\rho_{d}(Z,\bar Z;\eps))$ with
$\rho(p_0,\bar p_0;\eps)=0$ and
$\d\rho_1(\cdot;\eps)\wedge\cdots\wedge\d\rho_{d}(\cdot;\eps)\ne 0$ for every $\eps \in E$ 
such that there exists a neighbourhood of $\left\{ p_0 \right\} \times E\subset \CN\times E$
on which $\rho(Z,\bar Z; \eps)$ is real-analytic in all its arguments.
In particular, for each $\eps\in E$, the set $\{Z\in \C^N:\rho (Z,\bar Z;\eps)=0\}$ defines
a germ at $p_0$ of a real-analytic generic submanifold
$M_\eps\subset \C^N$ of codimension $d$. Given a fixed germ of a real-analytic generic submanifold
$M\subset \C^N$ through $p_0$ and $(M_\eps)_{\eps\in E}$  a real-analytic family through $p_0$ as defined above, 
we say that $(M_\eps)_{\eps\in E}$ is a real-analytic {\em deformation} of $(M,p_0)$ if
there exists $\eps_0\in E$ such that $(M,p_0)=(M_{\eps_0},p_0)$.

 We are now ready to state the
 following result.

\begin{Thm}\label{t:deformation}
Let $(M,p_{0})$ be a germ of a real-analytic generic submanifold of
codimension $d$ that is minimal and in the class ${\mathcal C}$ and
set $\ell_0 = 2 (d+1) \kappa_{M} (p_0)$. Let $(M_\eps)_{\eps \in E}$
be a real-analytic deformation of the germ $(M,p_0)$ $($parametrized by some real-analytic manifold $E)$  with $(M_{\eps_0},p_0)=(M,p_0)$
for some $\eps_0\in E$.
Then there exist open neighbourhoods $U_0$ of $\eps_0$ in $E$, 
$U_1$ of $p_0 \in \CN$ and $\Omega$ of $j_{p_0}^{\ell_0} {\sf Id}$ in $G^{\ell_0}_{p_0}(\C^N)$ and a real-analytic map 
$\Psi(Z, \L; \eps ) \colon U_1 \times \Omega \times U_{0} \to \CN$, 
holomorphic in its first factor such
that for every germ of a biholomorphic map $H\colon (\CN,p_0) \to
(\CN,p_0)$ sending $(M_\eps, p_0 )$ for some $\eps\in U_0$ into
$(M,p_0)$ with $j_{p_0}^{\ell_0} H \in \Omega$, we have
\[ H(Z) = \Psi (Z, j_{p_0}^{\ell_0} H; \eps ),\ {\rm for}\ Z\in \C^N\, {\rm close}\ {\rm to}\ p_0. \]
\end{Thm}

\begin{Rem}\label{cannot}
It is natural to ask whether Theorem~\ref{t:deformation} remains true
with the target manifold $(M,p_0)$ also varying. Such a result holds for 
finitely nondegenerate manifolds \cite{BER99b, KZ05}. 
However, it {\em cannot} hold for the more general class $\mathcal{C}$
(even in the real-analytic case) as the
example with $M\subset \C_{(z,w)}^2$ given by $\Im w=|z|^4$ shows,
see \cite[Example~1.5]{KZ05}.
\end{Rem}

Let us now show how Theorem~\ref{t:jetparamthm} follows from
Theorem~\ref{t:deformation}.

\begin{proof}[Proof of Theorem~{\rm \ref{t:jetparamthm}} assuming Theorem~{\rm \ref{t:deformation}}]
Without loss of generality, we may assume that $M$ is generic. Let
$\rho=\rho (Z,\bar Z)$ be a real-analytic vector valued defining
equation for $M$ in a neighbourhood $U$ of $0$. Consider the
real-analytic deformation of the germ $(M,p_0)$ obtained by varying
the base point in some small neighbourhood $\widetilde U\subset U$
i.e.\ defined by the family $(M_{p})_{p\in \widetilde U}$ where
$M_p$ is the germ at $p_0$ of the real-analytic generic submanifold
given by the equation $\{Z\in \C^N:\rho (Z-p_0+p,\bar
Z-\bar{p}_0+\bar p)=0\}$. Applying Theorem~\ref{t:deformation} to
this deformation, it is not difficult to derive the following:

\begin{Pro}\label{p:jetparaminv}
Under the assumptions of Theorem~{\rm \ref{t:jetparamthm}}, the set
of all germs $h\colon (M,p)\to (M,p_0)$ of local CR diffeomorphisms of
$M$ with variable source point $p$ has the real-analytic jet
parametrization property of order $\ell_0$ at $p_0$.
\end{Pro}

The conclusion of Theorem~\ref{t:jetparamthm} then follows easily
from Proposition~\ref{p:jetparaminv} and an application of the
inverse function theorem. We leave the details to the reader.
\end{proof}

\section{Proof of Theorem~\ref{t:deformation}} We assume that we are in the setting of
Theorem~\ref{t:deformation}. Without loss of generality, suppose
that $p_0$ coincides with the origin in $\C^N$ and set $n=N-d$.
Consider the given real-analytic 
family $(M_\eps)_{\eps\in E}$ and $\eps_0\in E$ satisfying
$(M_{\eps_0},0)=(M,0)$.

\subsection{Normal coordinates and Segre mappings for the deformation}

The first basic fact needed for the construction  of a mapping
 $\Psi$ satisfying the conclusion of Theorem~\ref{t:deformation}
 is the choice of a certain set of coordinates (the so-called
``normal coordinates'') for each manifold $M_\eps$ near the
origin and depending real-analytically on $\eps$ for $\eps$ close to $\eps_0$.

The coordinates we need are obtained from the standard construction of the
normal coordinates ({\em cf.} e.g. \cite{BERbook}):

\begin{Lem}
  \label{L:normalorabnormal} Let $(M_\eps)_{\eps\in E}$ be a real-analytic family of
  real-analytic generic submanifolds through the origin in $\C^N$
  of codimension $d$ and $\eps_0\in E$ as above. Then there exist
 germs of real-analytic maps
$$Z\colon (\CN\times E, (0,\eps_0))\to (\CN,0) \text{ and }
Q\colon (\Cn\times\Cn\times\Cd\times {E},(0,0,0,\eps_0))\to
(\Cd,0),$$ holomorphic in all their components except $E$,
 such that for every fixed $\eps\in E$ sufficiently close to $\eps_0$, the following holds:
  \begin{itemize}
    \item[(i)] ${Z}(0;\eps) = 0$ and the map $Z(\cdot;\eps)\colon (\CN,0)\to (\CN,0)$
is locally biholomorphic near $0$;
    \item[(ii)] in the local coordinates
      ${Z}(\cdot;\eps)=(z,w) \in \Cn\times \Cd$ near $0$,
      the manifold $M_\eps$ is given by
\begin{equation}\label{star}
      w - Q(z,\1{{z}}, \1{{ w}}; \eps) = 0;
\end{equation}
    \item[(iii)] one has $Q(z,0, \tau;\eps)
      \equiv Q(0, \chi, \tau;\eps) \equiv  \tau$.
  \end{itemize}
\end{Lem}

We fix an open neighbourhood $U_0$ of $\eps_0$ in $E$ so that
Lemma~\ref{L:normalorabnormal} holds. After possibly shrinking $U_0$
we may assume that for every $\eps\in U_0$, $M_\eps$ is minimal at
$0$ and that $\kappa_{M_\eps} (0) \leq \kappa_M (0)$; we set
$Q(z,\chi,\tau) := Q(z,\chi,\tau,\eps_0)$.

The next tools we need are the Segre mappings associated with
the manifolds $M_\eps$, $\eps \in U_0$. Recall that for every integer
$k\ge 1$, the $k$-th Segre (germ of a) mapping
$$v_\eps^k\colon (\C^{kn},0)\to (\C^n\times\Cd,0)$$
associated to $(M_\eps,0)$ and the chosen normal coordinates is
defined inductively as follows (see \cite{BER99b}):
\begin{equation}
  v_\eps^1 (t^{1}) := (t^1,0), \quad v_\eps^{k+1}(t^{[k+1]}) :=
  \left( t^{k+1}, Q\big(t^{k+1},\1{v_\eps^{k}} ( t^{[k]} ); \eps \big)\right),
  \label{e:segreparm}
\end{equation}
where $t^k \in \Cn$, $t^{[k]}: = (t^1,\dots,t^{k})\in\C^{kn}$. Here
and throughout the paper, for any power series mapping $\theta$, we
denote by $\bar{\theta}$ the power series obtained from $\theta$ by
taking complex conjugates of its coefficients.

For every $\eps \in U_0$ and a germ of a biholomorphic map $H\colon
(\CN,0) \to (\CN,0)$  sending $(M_\eps,0)$ into $(M,0)$, we define
\begin{equation}\label{e:given}
H_\eps \colon (\CN,0) \to (\CN,0), \quad H_\eps: = 
H(Z(\cdot;\eps)^{-1}),
\end{equation}
 where $Z(\cdot;\eps)^{-1}\colon (\CN,0)\to (\CN,0)$ is the local inverse of $Z(\cdot;\eps);$
 $H_\eps$ sends $(M_\eps,0)$ written  in the $Z$-coordinates into $(M,0)$, and
 $M_\eps$ is given by \eqref{star}. It is clear from the construction of the above
coordinates and from the Inverse Function Theorem that it is enough
to prove the parametrization property for all our mappings $H_\eps$
for $\eps$ sufficiently close to $\eps_0$ to obtain the conclusion of
Theorem~\ref{t:deformation}.

After choosing normal coordinates $Z'=(z',w')\in \C^n\times \C^d$
for the target manifold $M$ at $p_0$, which are fixed here, we write
$H_\eps=(f_\eps,g_\eps)\in \C^n\times \C^d$ and also denote by $\6M_{\eps}$ the germ at
the origin of the complexification of $M_\eps$. In our coordinates,
it is defined as the germ of the complex submanifold of $\C^N\times\C^N$ at the origin given by
$$\6M_\eps=\{(Z,\zeta)=((z,w),(\chi,\tau))\in \C^{n}\times \C^d\times
\C^{n}\times \C^d: w=Q(z,{ \chi},{\tau};\eps)\}.$$

\subsection{Reflection identities with parameters}\label{s:identities}

We now want to state a version with parameters of the reflection
identities given in \cite[Propositions 9.1 and 9.2, Lemma 9.3
and Proposition 9.4]{LM05}. For this, as in \cite{LM05},  it is
convenient to introduce the following notation.

For every positive integer $k$, we denote by $J_{0,0}^k(\C^N)$ the space of all jets at the
origin of order $k$ of holomorphic mappings from $\C^N$ into itself
and fixing the origin. In our normal
coordinates $Z = (Z_1,\dots,Z_N)$ in $\CN$, we identify a jet
${\mathcal J}\in J_{0,0}^{k} (\CN)$ with a polynomial map of the
form
\begin{equation}\label{jet}
{\mathcal J}={\mathcal J} (Z)=\sum_{\alpha \in \N^r,\
1\leq|\alpha|\leq k}\frac{\Lambda^k_{\alpha}}{\alpha !}\,
Z^{\alpha},
\end{equation}
where $\Lambda^k_{\alpha}\in \CN$. We thus have for a jet $
{\mathcal J}\in J_{0,0}^k(\C^N)$, the coordinates
\begin{equation}\label{e:jetnotation}
\Lambda^k:=(\Lambda^k_{\alpha})_{1\leq|\alpha|\leq k}
\end{equation} 
given by \eqref{jet}. Given a germ of a
holomorphic map $h\colon (\CN,0)\to (\CN,0)$, $h=h(t)$, for $t$
sufficiently small we use for the $k$-jet of
$h$ at $t$ the notation $j_t^kh=:(t,h(t),\widehat j_t^kh)$
(which is defined as a germ at $0$). Moreover, since
$h(0)=0$, we may also identify $j^k_0h$ with $\widehat j^k_0h$,
which we will freely do in the sequel.

Given the normal coordinates $(z,w)\in\Cn\times\Cd=\CN$, 
we consider a special component of a jet $\Lambda^k \in
J_{0,0}^k(\C^N)$ defined as follows. Denote the set of all multiindices of length one
having $0$ from the $n+1$-th to the $N$-th component by $S$, and the
projection onto the first $n$ coordinates  by ${\rm proj}_1\colon
\CN\to \Cn$ (that is, ${\rm proj}_1(z,w)=z$). Then set
\begin{equation}\label{e:crapagain}
\widetilde \Lambda^1:=({\rm proj}_1(\Lambda_{\alpha}))_{\alpha \in
S}.
\end{equation}
Note that for any local holomorphic map 
$$(\C^n \times \Cd,0) \ni
(z,w) \mapsto h(z,w)=(f(z,w),g(z,w))\in (\Cn \times \Cd,0),$$ if
$\jetm{0}{k}h=\Lambda^k$, then $\widetilde \Lambda^1=(\frac{\partial
f}{\partial z}(0))$. We can therefore identify $\widetilde
\Lambda^1$ with an $n\times n$ matrix or equivalently with an
element of $J_{0,0}^1(\Cn)$. Throughout the paper, given any jet
$\lambda^k\in J_{0,0}^k(\CN)$, $\widetilde \lambda^1$ will always
denote the component of $\lambda^k$ defined by \eqref{e:crapagain}.

In addition, for every positive integer $r$ and an open neighborhood $U_0$ of $\eps_0$ in $E$, 
we denote by $\6S_r=\6S_r(U_0)$ the ring of germs at $\{0\}\times U_0$ of 
real-analytic functions on $\C^r \times E$
that are holomorphic in their first argument.
Recall that this is
the space of all real-analytic 
functions that are defined in a connected open neighbourhood
(depending on the function) of $\{0\}\times U_0$ in $\C^r \times
E$ (and holomorphic in their first argument).

We now collect the following versions of the reflection identities of
\cite[Section 9]{LM05} with parameters that are necessary in order 
to complete the proof of Theorem~\ref{t:deformation}. 
The first basic identity given by  \eqref{e:fundamental1} 
is standard and may obtained by complexifying the identity $H_\eps(M_\eps)\subset M$ 
and applying the vector fields tangent to $\6M_\eps$.

\begin{Pro}\label{p:reflection1} In the above setting,
there exists a polynomial $\6D=\6D(Z,\zeta,\eps,
\Lambda^1)\in \6S_{2N}[\Lambda^1]$ and, for every $\alpha \in
\N^n\setminus \{0\}$, a $\Cd$-valued polynomial map
${\mathcal P}_{\alpha}={\mathcal
P}_{\alpha}(Z,\zeta,\eps,\Lambda^{|\alpha|})$ whose components are
in the ring $\6S_{2N}[\Lambda^{|\alpha|}]$ such that for $\eps \in
U_0$ and for every map $H_\eps\colon (M_\eps,0) \to (M,0)$ the
following holds:
\begin{enumerate}
\item[{\rm (i)}] $\6D (0,0,\eps,\Lambda^1)={\det}\, \widetilde
\Lambda^1$; 
\item[{\rm (ii)}] for all $(Z,\zeta)\in \M_\eps$
near $0$,
\begin{equation}\label{e:fundamental1}
(\6D(Z,\zeta,\eps,\widehat j_{\zeta}^{1} \overline{H}_\eps
))^{2|\alpha|-1}\,
\bar{Q}_{\chi^{\alpha}}(\bar{f}_\eps(\zeta),H_\eps(Z))={\mathcal
P}_{\alpha}(Z,\zeta,\eps,\widehat j_{\zeta}^{\left| \alpha \right|}
\overline{H}_\eps).
\end{equation}
\end{enumerate}
\end{Pro}

The next identity given by \eqref{e:fundamental2} involves the (transversal) 
derivatives of the mappings $H_\eps$
and follows easily from differentiating \eqref{e:fundamental1} and applying the chain rule.

\begin{Pro}\label{p:reflection2}
For any $\mu \in \N^d\setminus \{0\}$ and $\alpha \in
\N^n\setminus\{0\}$, there exist a $\Cd$-valued polynomial
map ${\mathcal
T}_{\mu,\alpha}(Z,\zeta,Z',\zeta',\eps,\lambda^{|\mu|-1},\Lambda^{|\mu|})$
whose components belong to the ring
$\6S_{4N}[\lambda^{|\mu|-1},\Lambda^{|\mu|}]$ and a
$\Cd$-valued polynomial map ${\mathcal Q}_{\mu,
\alpha}(Z,\zeta,\eps,\Lambda^{|\alpha|+|\mu|})$ whose components are
in the ring $\6S_{2N}[\Lambda^{|\alpha|+|\mu|}]$, such that for
$\eps \in U_0$, for every map $H_\eps\colon (M_{\eps},0)\to (M,0)$ 
and for any $(Z,\zeta)\in \6M_\eps$ close to the origin, the
following relation holds:
\begin{equation}\label{e:fundamental2}
\frac{\partial^{|\mu|}H_\eps}{\partial w^\mu}(Z)\cdot
\bar{Q}_{\chi^{\alpha},Z}(\bar{f}_\eps(\zeta),H_\eps(Z))
=(*)_1+(*)_2,
\end{equation}
where \begin{equation}\label{e:*1} (*)_1:= {\mathcal T}_{\mu,
\alpha}\left(Z,\zeta,H_\eps(Z),\overline{H}_\eps(\zeta),
\eps,\widehat j_{Z}^{\left| \mu \right| -1} H_\eps , \widehat
j_{\zeta}^{\left| \mu \right|} \overline{H}_\eps \right)
\end{equation} and
\begin{equation}\label{e:*2}
(*)_2:=\frac{{\mathcal Q}_{\mu,\alpha}(Z,\zeta, \eps, \widehat
j_{\zeta}^{|\alpha| + |\mu|} \overline{H}_\eps)}{(\6D(Z,\zeta,\eps,
\widehat j_{\zeta}^{1} \overline{H}_\eps ))^{2|\alpha|+|\mu|-1}}.
\end{equation}
\end{Pro}

In the next lemma, we observe that for any given map $H_\eps=(f_\eps,g_\eps)$, 
the (transversal) derivatives of the normal component $g_\eps$ can be expressed 
(in an universal way) through the (transversal) derivatives of the components of 
$f_\eps$ and some other terms that have to be seen as remainders. 
In particular, this lemma will allow us (as in \cite{LM05}) 
to derive the desired  parametrizations of the maps $H_\eps$ and their derivatives on each Segre set 
from the corresponding parametrizations of the maps $f_\eps$ and their derivatives.

\begin{Lem}\label{l:gderivative}
 For any $\mu \in \N^d\setminus \{0\}$, there exists a $\Cd$-valued polynomial map
 \[W_\mu={W}_{\mu}\left( Z,\zeta,Z',\zeta',\eps, \lambda^{|\mu|-1},\Lambda^{|\mu|} \right),\] 
whose components belong
to the ring $\6S_{4N}[\lambda^{|\mu|-1},\Lambda^{|\mu|}]$ and such that
for $\eps \in U_0$, for every map $H_\eps\colon (M_{\eps},0) \to
(M,0)$
and for any $(Z,\zeta)\in \6M_\eps$ close to the origin, the identity
\begin{equation}\label{e:gmu}
\frac{\partial^{|\mu|}g_\eps}{\partial
w^\mu}(Z)=\frac{\partial^{|\mu|}f_\eps}{\partial w^\mu}(Z)\cdot
Q_z(f_\eps(Z),\overline{H}_\eps(\zeta))+(*)_3
\end{equation}
holds with
\begin{equation}\label{e:*3}
(*)_3:=W_\mu \left(Z,\zeta,H_\eps(Z),\overline{H}_\eps(\zeta), \eps,
\widehat j_{Z}^{\left| \mu \right| - 1}H_\eps , \widehat
j_{\zeta}^{\left| \mu \right| } \overline{H}_\eps \right).
\end{equation}
\end{Lem}

The next statement is obtained as a direct combination of Lemma~\ref{l:gderivative} 
and Proposition~\ref{p:reflection2}
and provides the form of the system of equations fulfilled by any (transversal) derivative of $f_\eps$.

\begin{Pro}\label{p:reflection3}
For any $\mu \in \N^d\setminus \{0\}$ and $\alpha \in
\N^n\setminus\{0\}$, there exist a $\Cd$-valued polynomial
map ${\mathcal
T}'_{\mu,\alpha}(Z,\zeta,Z',\zeta',\eps,\lambda^{|\mu|-1},\Lambda^{|\mu|})$
whose components belong to the ring
$\6S_{4N}[\lambda^{|\mu|-1},\Lambda^{|\mu|}]$  such that for $\eps
\in U_0$ and for every map $H_\eps\colon (M_{\eps},0)\to(M,0)$  
the following relation holds for $(Z,\zeta)\in \6M_\eps$ close to
$0$:
\begin{equation}\label{e:fundamental3}
\frac{\partial^{|\mu|}f_\eps}{\partial w^\mu}(Z)\cdot
\left(\bar{Q}_{\chi^{\alpha},z}(\bar{f}_\eps(\zeta),H_\eps(Z))+
Q_z(f_\eps(Z),\overline{H}_\eps(\zeta))\cdot
\bar{Q}_{\chi^{\alpha},w}(\bar{f}_\eps(\zeta),H_\eps(Z))\right)=(*)_1'+(*)_2,
\end{equation}
where $(*)_2$ is given by \eqref{e:*2} and $(*)_1'$ is given by
\begin{equation}\label{e:*1'}
(*)_1':={\mathcal
T}_{\mu,\alpha}'\left(Z,\zeta,H_\eps(Z),\overline{H}_\eps(\zeta),\eps,
\widehat j_{Z}^{\left| \mu \right| - 1} H_\eps, \widehat
j_{\zeta}^{\left| \mu \right|}{\overline{H}_\eps} \right).
\end{equation}
\end{Pro}

Since the proof of the above relations is analogous to those derived in
\cite{LM05}, we leave the details to the reader.  We should point out that, in the reflection identities with parameters mentioned above, the most relevant fact is the location of the parameter $\eps$ in the identities. Indeed, the parameter $\eps$ appears always in an appropriate place so that the results concerning the parametrization of solutions of singular analytic systems given in the next paragraph will be applicable. This crucial fact explains why we can follow
the analysis of \cite{LM05} in order to derive Theorem~\ref{t:deformation}.

\subsection{Parametrization of solutions of singular analytic systems}\label{s:parametrization}

We state here the two versions of the parametrization results for
singular systems needed for the proof of Theorem~\ref{t:deformation}.
The first one is needed to have a parametrization of the compositions
$H_\eps\circ v^j_\eps$ for all integers $j$,
where $v^j_\eps$ is defined by \eqref{e:segreparm}.  

\begin{Thm}
  \label{c:corparcn2}   Let $A:(\C^m,0)\to \C^m$ be
  a germ of a holomorphic map
 of generic rank $m$, $X$ a real-analytic
 manifold, $Y$ a complex manifold
 and $b=b(z,x,y)$ a $\C^m$-valued real-analytic map defined
  on an open neighbourhood $V$ of $\{0\}\times X\times Y$ in $\C^m\times X\times Y$,
holomorphic in $(z,y)$. 
Then there exists a 
real-analytic map $\Gamma=\Gamma(z,\genmat,x,y)\colon
\C^m\times\glmc\times {X}\times Y\to \C^m$, defined on an open
neighbourhood $\Omega$ of $\{0\}\times \glmc\times {X}\times Y$,
holomorphic in all its components except $X$, satisfying the
following properties:
\begin{enumerate}
\item[(i)]  If $u:(\C^m,0)\to(\C^m,0)$ is a germ of a
biholomorphism satisfying
  $A(u(z)) = b(z,x_{0},y_0)$ for some $(x_{0},y_0)\in {X}\times Y$, then necessarily
  $u(z) = \Gamma(z,j_0^1u,x_0,y_0)$;
  \item[(ii)] For every $\lambda \in \glmc $ and $(x_0,y_0)\in X\times Y$, the map $\Gamma$
  satisfies $\Gamma (0,\lambda,x_0,y_0)=0$ and $\displaystyle \frac{\partial \Gamma}{\partial z}
  (0,\lambda,x_0,y_0)=\lambda$.
  \end{enumerate}
\end{Thm}

The statement given by Theorem~\ref{c:corparcn2} follows directly from \cite[Corollary 3.2]{LM05} after an obvious  complexification argument.

The second version given below is needed to get a parametrization of the mappings
$(\partial^\beta H_\eps)\circ v^j_\eps$ for all integers $j$ and all
multiindices $\beta \in \N^N \setminus \{0\}$. 

\begin{Pro}
  \label{t:lineqn2}
  Let $\Theta$ be an $r\times r$ matrix with
  holomorphic coefficients near the origin in $\C^m$, $m,r\geq 1$,
  such that $\Theta$ is of generic rank $r$. Let $X$ be a real-analytic manifold and $Y$
  a complex manifold. Assume that
  $c\colon \C^m\times X\times Y \to \C^m$ and $b\colon \C^m\times X\times Y \to\C^r$
  are real-analytic maps defined on some
  neighbourhood $V$ of $\{ 0 \}\times X\times Y$ such that
   $(z,y)\mapsto b(z,x,y)$ and $(z,y)\mapsto c(z,x,y)$
    are holomorphic on $V_{x} = \{ (z,y)\in\C^m \times Y  \colon (z,x,y) \in V \}$ for each $x\in
   X$. Assume furthermore that $c$ satisfies
  \[ c(0,x,y) = 0, \quad \det \, c_z (0,x,y) \neq 0, \quad {\rm for}\ {\rm every}\ (x,y)\in X\times Y. \]
  Then there exists a real-analytic map
  $\Gamma\colon \C^m\times X\times Y\to \C^r$ defined on a neighbourhood of $\{0\} \times
  X\times Y$, holomorphic in all its components except $X$,
  such that if $u\colon (\C^m,0)\to\C^r$
  is a germ of a holomorphic map satisfying $ \Theta(c(z,x_0,y_0)) \cdot u(z) = b(z,x_0,y_0)$
  for some $(x_0,y_0)\in X\times Y$, then
  $u(z)= \Gamma(z,x_{0},y_0)$.
\end{Pro}

The statement given by Proposition~\ref{t:lineqn2} follows from 
\cite[Proposition 6.3]{LM05} and again a simple complexification argument.

\subsection{Completion of the proof of Theorem~\ref{t:deformation}}
With the statements given in Sections
\ref{s:identities} and \ref{s:parametrization} at our disposal, we
can follow the plan of the proof of \cite[Theorem 7.3]{LM05} to get
the needed parametrization of the maps $H_\eps$ restricted to any
Segre set. More precisely, the reader may verify  that after applying 
the above statements as in
\cite{LM05}, one obtains the following:

\begin{Pro}
In the above setting and shrinking the neighbourhood $U_0$ if
necessary, for every positive integer $j$, there exists a real-analytic map
\[ \Psi_j \colon \C^{nj}\times {U_0} \times
J_{0,0}^{j \kappa_M (0)} (\CN)\to \CN,
\] defined in a neighbourhood of $\{ 0 \} \times U_0 \times W_j
$ where $W_j$ is an open set in the jet space containing all the
jets $($at $0)$ of the maps $H_\eps$ for $\eps \in U_0$, that is
holomorphic in its first factor and satisfying in addition
\begin{equation}\label{e:sicksicksick}
  \left( H_\eps \circ v_\eps^{j}\right) \left(  t^{ [j]} \right) =
\Psi_j \left( t^{[j]}, \eps, j_0^{j \kappa_M (0)} H_\eps\, \right),
\end{equation}
for all $t^{[j]}$ sufficiently close to the origin.
\end{Pro}

We are now in a position to finish the proof of
Theorem~\ref{t:deformation}. 
For this, recall first that $\ell_0=2(d+1)\kappa_M(0)$ and consider the equation
\eqref{e:sicksicksick} for $j = 2 (d +1)$ that we localize near the
point $(\eps_0,j_0^{\ell_0}{\sf Id})\in E \times J_{0,0}^{\ell_0} (\CN)$. 
Shrinking $U_0$ if necessary, there exist open
neighbourhoods $O\subset \C^{2n(d+1)}$ of the origin and $O'\subset
J_{0,0}^{\ell_0} (\CN)$ of $j_0^{\ell_0}{\sf Id}$ such
that $\Psi_{2(d+1)}$ is defined over $O\times U_0\times O'$ and such
that  for every $\eps \in U_0$ satisfying $j_0^{\ell_0}H_\eps \in
O'$, the identity \eqref{e:sicksicksick} holds (with $j=2(d+1)$) for
$t^{[2(d+1)]}$ sufficiently close to the origin.

The rest of the proof closely follows the lines of  \cite[Section~4]{KZ05}; it
 consists of using a version of the implicit function with singularities \cite[Lemma~3.4]{KZ05}
and resolving the obtained singularities by using \cite[Lemma~4.3]{KZ05}. The differences between 
the situation treated in the present paper and that of \cite{KZ05} are the parameter dependence which is real-analytic in our case (instead of smooth in \cite{KZ05}) and the absence of the error terms 
in the formula \eqref{e:sicksicksick} (in contrast to \cite{KZ05}). 
The details are left to the reader. 

\section{Proof of Theorem~\ref{t:unify}}\label{ss:final}

\begin{proof}[Proof of Theorem~{\rm \ref{t:unify}}]
Suppose first that $M$ is a connected real-analytic CR-submanifold
in $\C^N$. Then we claim that the conclusion of
Theorem~\ref{t:unify} follows from the conjunction of
Theorem~\ref{main-technical} and Theorem~\ref{t:jetparamthm}.
Indeed, assumption (i) of Theorem~\ref{t:unify} and
Theorem~\ref{t:jetparamthm} imply that assumption (i) of
Theorem~\ref{main-technical} is satisfied. (Note that the upper semi-continuity 
of the integer $\kappa_M(p)$ on $p\in K\subset M$ in Theorem~\ref{t:jetparamthm} is also used here in order to deduce the existence of the integer $k$ satisfying the conclusions of Theorem~\ref{main-technical} (i)). 
Furthermore, assumption
(ii) of Theorem~\ref{t:unify} together with the results of
\cite{BER99b, KZ05} imply that assumption (ii) of
Theorem~\ref{main-technical} is also satisfied. This proves the
claim.

If $M$ is not connected, we may repeat the arguments of the proof of
\cite[Theorem 6.2]{BRWZ} since $M$ is assumed to have finitely many
connected components.

Finally, when $M$ is an abstract real-analytic CR manifold, the proof is
the same as before since it is based on purely local arguments and since any such manifold can locally be embedded 
as a CR submanifold of some complex euclidean space $\C^q$  for some integer $q$ (see e.g.\ \cite{BERbook}). The
proof of the theorem is complete.
\end{proof}

\section{An elementary proof of
Corollary~\ref{t:liehighcodim}}\label{s:lmz}

We conclude this paper by providing an elementary proof of
Corollary~\ref{t:liehighcodim} which avoids the use of
Theorem~\ref{main-technical} and rather follows the proof of
\cite[Corollary 1.3]{Z97}. Note that in any case one has to make use
of Theorem~\ref{t:jetparamthm}.

\begin{proof}[Proof of Corollary~{\rm \ref{t:liehighcodim}}]
Since $M$ is compact (and everywhere minimal) and embedded in some
Stein manifold, we may apply Theorem~\ref{t:jetparamthm} to conclude
that there exists a finite number of points $p_1,\ldots,p_k\in M$
and open neighbourhoods $\Omega'_j$ of $p_j$ in $M$ covering $M$
such that for every $h\in\Aut_\CR(M)$ sufficiently close to the
identity mapping, say in an open neighbourhood $\mathcal{N}$ of it,
Theorem~\ref{t:jetparamthm} holds at all points $p_j$ with a
parametrization $\Psi_j$
  defined in a neighborhood of $\Omega'_j\times\{p_j\}$ with the jet order
  $\ell_j$. Write $\ell =\max \ell_j$.
 As in
  \cite{Z97}, our goal is to show that the
  image of the neighbourhood $\mathcal{N}\subset \Aut_\CR(M)$ under the
  homeomorphism (onto its image)
  \[ h \mapsto \eta(h) = \left(
  j_{p_1}^{\ell} h, \ldots, j_{p_k}^{\ell} h,
  j_{p_1}^{\ell} h^{-1},\ldots, j_{p_k}^{\ell} h^{-1}\right) \in
    \left( G_{p_1}^{\ell} (M)\times \ldots \times
    G_{p_k}^{\ell} (M)\right)^2 =: \mathcal{Y}^2
  \]
  is a real-analytic subset of the target space, and that the
  group law is real-analytic. But it is
  easy to single out the points in the image which give rise to
  a global automorphism of $M$. Any
  $ (\alpha,\beta)=
  \left( \alpha_1,\ldots, \alpha_k, \beta_1,\ldots,\beta_k
  \right) \in \mathcal{Y}^2$  belongs to  $\eta ({\mathcal N})$
  if and only if for every $j, m = 1,\ldots, k$,
  the following identities are satisfied:
  \[\begin{gathered}
    \Psi_j( \cdot , \alpha_j) = \Psi_m (\cdot, \alpha_m), \,
    \Psi_j( \cdot , \beta_j) = \Psi_m (\cdot, \beta_m)
    \text{ on }  \Omega'_j \cap \Omega'_m, \\
    \Psi_m \left( \Psi_{m}
    \left( \cdot, \alpha_m \right),\beta_m \right) =
       \Psi_m \left( \Psi_{m}
       \left( \cdot, \beta_m \right),\alpha_m \right) = {\sf Id} \quad  \text{ near } p_m, \\
     \alpha_j = j_{p_j}^{\ell} \left(\Psi_j (\cdot, \alpha_j) \right), \,
     \beta_m = j_{p_m}^{\ell} \left( \Psi_m (\cdot, \beta_m)
     \right).
  \end{gathered}\]
  From this, it is clear that $\eta ({\mathcal N})$ is a real-analytic
  subset of $\mathcal{Y}^2$, and again following~\cite{Z97}, we
  see that the group law is indeed real-analytic.
  This concludes the proof of
  Corollary~\ref{t:liehighcodim}.
\end{proof}

\end{document}